\magnification=\magstep1
\documentstyle{amsppt}
\hsize=5.5truein
\vsize=7.5truein
\NoBlackBoxes
\NoPageNumbers

%\input amstex
%\documentstyle{amsppt}
%\magnification=1200
%\hsize=5.5trueinch
%\vsize=7.5trueinch
\def\deg{\text{deg}\thinspace}
\def\Proj{\text{Proj}\thinspace}

%-THIS VERSION WAS WRITTEN ON SEPTEMBER 17TH 1999, WHEN I WANTED TO ADD THE --
%------EGA LEMME TO THE PAPER. THE TEX FILE LOOKS A BIT DIFFERENT FROM THE----
%------PREVIOUS ONES----------------------------------------------------------

%---------------------------------------------------------------------------  
\topmatter

\title Normal Ideals of Graded Rings \endtitle
\author Sara Faridi \endauthor
\dedicatory Mathematics Department, University of Michigan, Ann Arbor, MI 48109 \\email: faridi\@math.lsa.umich.edu
\enddedicatory
%\address Mathematics Department, University of Michigan, Ann Arbor, MI 48109 
%\endaddress
%\email faridi\@math.lsa.umich.edu \endemail
%\date January 1999 \enddate
\abstract For a graded domain $R=k[X_0,...,X_m]/J$ over an arbitrary domain $k$, it is shown that the ideal generated by elements of degree 
$\geq mA$, where $A$ is the least common multiple of the weights of the $X_i$, is a normal ideal.  \endabstract 

\thanks I would like to thank my advisor Professor Karen Smith for all her help and many fruitful discussions during 
the preparation of this paper. \endthanks
\endtopmatter

%----------------------------------------------------------------------------
\document 
 
\settabs 4 \columns
  
\head 1. Introduction \endhead

  	In the theory of resolution of singularities, one wishes to be
 able to blow up a singular variety along a closed subscheme and
 obtain a smooth variety birational to the original one.  One question
 that comes up is the existence of such resolutions; it is known for
 varieties over fields of characteristic zero [H], and is conjectured
 in the prime characteristic case.  Another question is describing the
 closed subschemes that give smooth blowups.  Translating this problem
 into the language of algebra, one is interested in ideals $I$ of a
 ring $R$ such that $\Proj R[It]$ is a smooth variety, where $R[It]=
 \oplus_{n \in \bold N} {I^nt^n}$ is the Rees ring of $R$ along $I$.
 This leads one to focus on the Rees ring and its properties.  In
 particular one would like to know when the Rees ring is normal. If
 $R$ is a normal domain, $R[It]$ is normal if and only if $I$ is a
 normal ideal, where by $ normal$ $ideal$ we mean an ideal all of
 whose positive powers are integrally closed.  Normal ideals have been
 studied in different cases, and in some special cases necessary and
 sufficient conditions have been given for an ideal to be normal; see
 for example [G], [O], and [HS].

  	In general given a ring, finding an ideal whose blowup is
regular has proven to be a difficult problem.  The same is true even
for finding normal ideals.  In this paper we will construct a normal
ideal $I$ for a general graded domain $R=k[X_0,...,X_m]/J$ which
depends only on the weights of the variables $X_0,...,X_m$, and is
thus very simple to construct. We will begin by recalling some
necessary definitions. For more on the construction of Rees rings see
[E].

	Let $R$ be any $\bold N$-graded domain, which is a quotient of
 a polynomial ring $k[X_0,...,X_m]$ modulo a homogeneous ideal $J$,
 where $k=R_0$ is an arbitrary domain and $X_0,...,X_m$ are variables
 of positive weights $A_0,...,A_m$.  In practice, $k$ is usually a
 field. Let $\bold m =(X_0,....,X_m)$ be the irrelevant ideal of
 $R$. By abuse of notation, by $X_0,...,X_m$ we will mean the images
 of the variables $X_0,...,X_m$ in $R$.  Throughout this paper
 $R_{\geq \alpha}$ refers to the ideal of $R$ generated by the
 elements of degree at least $\alpha$ in the graded ring $R$.

	We begin our search for normal ideals by reviewing some basic
	facts about them.

   \definition {Definition 1}
	For Noetherian domains $S \subseteq T$ the $integral$
	$closure$ of $S$ in $T$, denoted by $\overline S$, is defined
	as all $x \in T$ that satisfy an equation of the form
	$x^n+a_1x^{n-1}+...+a_{n-1}x+a_n=0$ where $a_i \in S$, for
	$i=1,...,n$.
\enddefinition

   \definition {Definition 2} For a Noetherian domain $S$ and an ideal
	$J$ of $S$, the $integral$ $closure$ of $J$ in $S$, denoted by
	$\overline J$, is defined as

	 all $x \in S$ that satisfy an
	equation of the form $x^n+a_1x^{n-1}+...+a_{n-1}x+a_n=0$ where
	$a_i \in J^i$, $i=1,...,n$; or equivalently,

  	all $x \in S$ for which there exists a nonzero $c \in S$ such
  	that $cx^n \in J^n$ for all positive integers $n$ (see [Ho]).
\enddefinition

   \demo{Note} These definitions also hold in the nondomain case; we
   just have to choose $c$ not in any minimal prime of $S$ in the
   second statement. \enddemo
  
        Here are a few well-known facts that we shall use:
 
\proclaim {Theorem 1} If S is a normal domain, then $S[It]$ is normal iff
 $I^n$ is integrally closed for every positive integer $n$.\endproclaim

\proclaim {Theorem 2} In an $\bold N$-graded domain $(S,\bold n)$, for
 any positive integer $\alpha$, the $\bold n$-primary ideal $I=S_{\geq
 \alpha}$ is integrally closed.
\endproclaim

           Using Theorems 1 and 2, we will be looking for an ($\bold m
        $-primary) ideal $I$ of the form $R_{\geq \alpha}$, with the
        property that $I^n=R_{\geq n\alpha}$ for all integers $n\geq
        1$. Note that from Lemma (2.1.6) in [EGA] one can deduce that
        such an $\alpha$ always exists. Theorem 3 gives an explicit
        value for $\alpha$.

%----------------MAIN THEOREM: THEOREM 3---------------------------------

\head 2. Main Theorem \endhead

%\vskip .2cm
 
\proclaim {Theorem 3} Let $R$ be a graded domain, which is a quotient of a 
 polynomial ring $k[X_0,...,X_m]$ modulo a homogeneous ideal $J$,
 where $k$ is an arbitrary domain and $X_0,...,X_m$ are variables of
 positive weights $A_0,...,A_m$. Let $A$ be the least common multiple
 of $A_0,...,A_m$. Then the ideal $I=R_{\geq mA}$ is a normal ideal.
 In particular, if $R$ is normal, the Rees ring $R[It]$ is
 normal.\endproclaim

  	We will show that for all positive integers $p$,
 $I^p=R_{\geq pmA}$. By Theorem 2, this will complete the proof.

%-----------------INDUCTIVE STEP------------------------------------------

   \proclaim {Inductive Step} Let $R$, $I$, and $A$ be as in Theorem
  3.  For $p \geq 2$, if $I^{p-1}=R_{\geq (p-1)mA}$ then $I^p=R_{\geq
  pmA}$.
\endproclaim

	It is obvious that $I^p \subseteq R_{\geq pmA}$.  We have to
	show that the other inclusion holds.  Also, observe that
	elements in $R_{\geq pmA}$ are sums of monomials in the $X_i$
	with coefficients in $k$, whose degrees are larger than or
	equal to $pmA$, and therefore we notice that $R_{\geq pmA}$ is
	generated by such monomials.  Hence it is enough to show that
	every monomial of $k[X_0,...,X_m]$ which lies in $R_{\geq
	pmA}$ also lies in $I^p$.
 
%---------------LEMMA AND PROOF--------------------------------------------
 
 \proclaim{Lemma} With notation as above, let $X_0^{c_0}...X_m^{c_m}
 \in R_{\geq pmA}$, where $c_0,...,c_m$ are nonnegative integers, and
 let $a_0,...,a_m$ be positive integers such that $a_iA_i=A$ for all
 $i$. Assume $I^{p-1}=R_{\geq (p-1)mA}$. Fix $n$ such that $1 \le n
 \le m$, and suppose that $a_i \le c_i$ for $0 \le i <n$. For each $i$
 smaller than $n$, let $k_i$ be the unique positive integer such that
 $k_ia_i \le c_i < (k_i+1)a_i$.  Then:

   1) if $k_0+...+k_{n-1} \le m-1$, then for some $a_j$ $(n \le j \le
   m)$, $a_j \le c_j$;
  
   2) if $k_0+...+k_{n-1} \geq m$, then $X_0^{c_0}...X_m^{c_m} \in I^p$.

\endproclaim
  
 \demo{Proof of Lemma}	$1)$	 If $c_i<a_i$ for all $n \le i \le m$, then 
$$
\alignat 1
 pmA &\le c_0A_0+...+c_mA_m \\
&< (k_0+1)a_0A_0+...+(k_{n-1}+1)a_{n-1}A_{n-1}+a_nA_n+...+a_mA_m \\        
&=(k_0+...+k_{n-1})A +nA+(m-n+1)A \\
&\le (m-1)A+nA+(m-n+1)A \\
&= 2mA.
\endalignat
$$

	It follows that $p<2$, contrary to the hypothesis.

   $2)$	Choose nonnegative integers $s_0,...,s_{n-1}$ such that 
$s_i \le k_i$ for all $0 \le i \le n-1$, and $s_0+...+s_{n-1}=m$.  We see 
that                                                                  
 $$X_0^{c_0}...X_m^{c_m} = X_0^{s_0a_0}...X_{n-1}^{s_{n-1}a_{n-1}} \cdot
 X_0^{c_0-s_0a_0}...X_{n-1}^{c_{n-1}-s_{n-1}a_{n-1}}X_n^{c_n}...X_m^{c_m}.$$  
	Now 
$$
\alignat 1
&\deg   X_0^{s_0a_0}...X_{n-1}^{s_{n-1}a_{n-1}} \\ 
&=s_0A+...+s_{n-1}A \\
&=mA,
\endalignat
$$                                                                             which implies that  $X_0^{s_0a_0}...X_{n-1}^{s_{n-1}a_{n-1}} \in I$.                                    
   	On the other hand 
$$
\alignat 1
 &\deg   X_0^{c_0-s_0a_0}...X_{n-1}^{c_{n-1}-s_{n-1}a_{n-1}}X_n^{c_n}
...X_m^{c_m} \\
&= \deg X_0^{c_0}...X_m^{c_m}-mA \\
&\geq pmA - mA  \\
&=(p-1)mA. 
\endalignat
$$

  	Therefore $ X_0^{c_0-s_0a_0}...X_{n-1}^{c_{n-1}-s_{n-1}a_{n-1}}
X_n^{c_n}...X_m^{c_m} \in I^{p-1}$, so $X_0^{c_0}...X_m^{c_m} \in I^p$. 
 \hfill $\blacksquare$

\enddemo  
                                                                    
%----------------PROOF OF INDUCTIVE STEP-------------------------------------  

\demo{Proof of Inductive Step}	We have $X_0^{c_0}...X_m^{c_m} \in 
R_{\geq pmA}$, where $c_0,...,c_m$ are nonnegative integers, and we
want to show that it lies in $I^p$ as well. Let $a_0,...,a_m$ be
positive integers such that $a_iA_i=A$ for all $i$.

  	If $c_i \geq ma_i$ for any $i$, say if $c_0 \geq ma_0$, then
	we will get $$X_0^{c_0}...X_m^{c_m} = X_0^{ma_0} \cdot
	X_0^{c_0-ma_0}...X_m^{c_m}.$$ Now, $\deg
	X_0^{ma_0}=ma_0A_0=mA$, and so $X_0^{ma_0} \in I$.  On the
	other hand,

 $$\deg X_0^{c_0-ma_0}...X_m^{c_m}= \deg X_0^{c_0}...X_m^{c_m}-mA
\geq pmA-mA=(p-1)mA.$$
	
	So $X_0^{c_0-ma_0}...X_m^{c_m} \in I^{p-1}$ by our
 assumption.  Therefore $X_0^{c_0}...X_m^{c_m} \in I^p,$ and we are done.                                                          

  	So let us look at the case where $c_i < ma_i$ for all $i$.

  	Now, we cannot have $c_i < a_i$ for all $i$, because in that
  	case we get
  
$$pmA \le c_0A_0+...+c_mA_m < a_0A_0+...+a_mA_m=(m+1)A,$$                     
and so  $$pm<m+1 \Rightarrow p<1+{1 \over m} \le 2,$$          		
which is a contradiction.                                                                     

  	So without loss of generality let $c_0 \geq a_0$, and let $k_0$ be the positive integer ($0<k_0<m$) such that $k_0a_0 \le c_0 < (k_0+1)a_0$. 

  	Hereafter, we proceed by induction: Having chosen
 $k_0,...,k_{n-1}$, if $k_0+...+k_{n-1} \geq m$, then
 $X_0^{c_0}...X_m^{c_m} \in I^p$ by the above lemma, and hence we are
 done.  

	However, if $k_0+...+k_{n-1} \le m-1$, by part 1 of the lemma
 there exists a $j \geq n$ for which $a_j \le c_j$. We can without
 loss of generality assume that $j=n$, and repeat the same cycle. Thus
 we are reduced to considering the case where $a_i \le c_i$ for
 $i=0,...,m$.  Let the $k_0,...,k_m$ be the unique positive integers
 for which $k_ia_i \le c_i < (k_i + 1)a_i$, $(0 \le i \le m)$. We can
 easily see that $k_0+...+k_m \geq m$, because otherwise
$$
\alignat 1
pmA &\le c_0A_0+...+c_mA_m \\ &< (k_0+1)a_0A_0+...(k_m+1)a_mA_m \\ &
=(k_0+...+k_m)A +(m+1)A \\ & \le (m-1)A+(m+1)A \\ &= 2mA,
\endalignat
$$    
which implies that $p <2$, contrary to our assumption. It follows then
from part 2 of the same lemma that $X_0^{c_0}...X_m^{c_m} \in I^p$.
This completes the proof of the inductive step.  \hfill $\blacksquare$
\enddemo                                                                   
%--------------------PROOF OF THEOREM 3------------------------------------                                                                   
   \demo{ Conclusion of Proof of Theorem 3} We need $I^p$ to be
   integrally closed for all $p \geq 1$. By Theorem 2, it suffices to
   have $I^p=R_{\geq pmA}$ for $p \geq 1$. We prove this by
   induction. The case $p=1$ is the definition of $I$. If $I^k=R_{\geq
   kmA}$ for all $1 \le k \le p-1$, then from the inductive step it
   follows that $I^p=R_{\geq pmA}$.  If $R$ is normal, it follows from
   Theorem 1 that $R[It]$ is a normal ring.\hfill $\blacksquare$
\enddemo

%-------------------EXAMPLES AND REMARKS-----------------------------------
\head 3. Examples \endhead

\example {Example 1} Let $k$ be a field, and $R=k[x,y,z]/(x^2+y^3-z^5)$,
 where $x,y,z$ have weights 15,10, and 6, respectively. the least
 common multiple of these variables is 30, and therefore by theorem 3,
 $$ \split I&=R_{\geq 60}= (
 x^4,x^3y^2,x^3yz,x^3z^3,x^2y^3,x^2y^2z^2,x^2yz^4,x^2z^5,xy^5, \\
 &xy^4z,
 xy^3z^3,xy^2z^5,xyz^6,xz^8,z^{10},yz^9,y^2z^7,y^3z^5,y^4z^4,y^5z^2,y^6)
 \endsplit $$ is a normal ideal for this ring. \endexample
\example {Example 2} In general, if $k$ is a field, and
 $R=k[x,y,z]/(x^a+y^b+z^c)$ is a domain, where the variables $x, y,z$
 have weights $bc,ac,$ and $ab$, respectively, the ideal $I=R_{\geq
 2abc}$ will be a normal ideal of $R$. \endexample

\example {Example 3} For the polynomial ring $R=k[x,y,z]$ over a field $k$,
 one can find many normal ideals by assigning different weights to the
 variables $x$, $y$ and $z$. For example, if we set $\deg x=1$, $\deg
 y=1$ and $\deg z=2$ we find that the ideal $I=R_{\geq
 4}=(x^4,x^3y,x^2y^2,x^2z,xy^3,xyz,y^4,y^2z,z^2)$ is a normal
 ideal. \endexample

\example {Remark 1} The method described in Example 3 above for
 finding normal ideals of polynomial rings is not as interesting in
 the case of two variables. It is a theorem due to Zariski (see [Z]),
 that the set of integrally closed ideals in a regular local ring of
 dimension two is closed under multiplication. Therefore, in the ring
 $R=k[x,y]$, all integrally closed ideals are normal. In particular,
 all ideals of the form $I=R_{\geq \alpha}$, where $\alpha$ is a
 positive integer, are normal (all powers of $I$ are integrally
 closed). \endexample

\example {Remark 2} Looking at Theorems 1 and 2, it is natural to wonder
 whether for an $\bold N$-graded normal ring $R$ and all positive
 integers $\alpha$, one can say that $\overline {(R_{\geq \alpha})^n}
 = R_{\geq n \alpha}$ for $n \in \bold N$. The answer in general is
 no. For example, consider the ring $R=k[x,y]$ where $k$ is a field,
 $\deg x=2$, and $\deg y =3$. Let $I={(R_{\geq 7})^3}$. We will check
 if $\overline I$ is equal to $R_{\geq 21}$. Since $R_{\geq 7}=
 (x^4,yx^2,xy^2,y^3)$, we can calculate
 $I=(x^{12},x^{10}y,x^8y^2,x^6y^3,x^5y^4,x^4y^5,x^3y^6,x^2y^7,xy^8,y^9)$. If
 $\Gamma$ is the set of pairs $(a,b)$ corresponding to generators
 $x^ay^b$ of $I$, we plot the points in $\Gamma$ on the $\bold R^2$
 plane, and we see that there are no pairs of positive integers
 $(c,d)$ in the convex hull of the region $\Gamma + {\bold R^2}^+$,
 besides those in $\Gamma + {\bold R^2}^+$ itself. This implies that
 $I$ is an integrally closed ideal (see [E]). On the other hand
 $x^{11} \in R_{\geq 21}$, but $x^{11}$ does not belong to $ I=
 \overline I$. Therefore $R_{\geq 21}$ cannot be the integral closure
 of $I$. (Note that it also follows from Zariski's theorem mentioned
 in Remark 1 that $I$ is an integrally closed ideal.) \endexample

\example {Remark 3} One might ask if all ideals of the form   
$I=R_{\geq \alpha}$, $\alpha \in {\bold N}$, in a normal graded ring
are normal.  The answer is no. A counterexample is the ring
$R=k[x,y,z]/(x^2+y^3z+z^4)$, where $x$, $y$ and $z$ have degrees 2, 1
and 1, respectively. The ideal $I=R_{\geq 1} =(x,y,z)$ is integrally
closed, but $I^2= (x^2,xy,xz,y^2,yz,z^2)$ is not:
$(x)^2+(y^2)(yz)+(z^2)^2=0$, and hence $x \in \overline {I^2}$, but
$x$ is not in $I^2$.  This, by the way, is an example of a ring whose
test ideal is not normal, given by Hara and Smith in [HSm]. For more
on test ideals and tight closure theory in general, see [HH].
\endexample

%-----------------------REFERENCES----------------------------------------

\enddocument